\UseRawInputEncoding
\documentclass[oneside, 11pt]{amsart}
\usepackage{amsmath}
\usepackage{amssymb}
\usepackage{dsfont}
\usepackage{color}
\usepackage[usenames,dvipsnames,x11names,svgnames]{xcolor}
\usepackage[colorlinks=Black,linkcolor=Black,citecolor=Black, urlcolor=Blue]{hyperref}
\usepackage{amsthm}
\usepackage{amscd}
\usepackage{geometry}
\usepackage{mathrsfs}
\usepackage{mathtools}

\newtheorem{theorem}{Theorem}

\newtheorem{proposition}{Proposition}

\newtheorem{corollary}{Corollary}

\theoremstyle{definition}

\newtheorem*{remark}{Remark}

\newcommand{\be}{\begin{equation*}}
\newcommand{\ee}{\end{equation*} }

\newcommand{\ben}{\begin{equation}}
\newcommand{\een}{\end{equation} }

\newcommand{\bs}{\begin{split}}
\newcommand{\es}{\end{split}}

\newcommand{\bmu}{\begin{multline*}}
\newcommand{\emu}{\end{multline*}}

\newcommand{\bmun}{\begin{multline}}
\newcommand{\emun}{\end{multline}}

\begin{document}



\title[]{On the Riemann-Hardy hypothesis for the Ramanujan zeta function}

\author{Xiao-Jun Yang$^{1,2}$}

\email{dyangxiaojun@163.com; xjyang@cumt.edu.cn}

\address{$^{1}$ School of Mathematics,and State Key Laboratory for Geo-Mechanics and
Deep Underground Engineering, China University of Mining and Technology, Xuzhou 221116, China}
\address{$^{2}$ Department of Mathematics, Faculty of Science, King Abdulaziz University P.O. Box 80257, Jeddah 21589, Saudi Arabia}

\begin{abstract}
The Ramanujan zeta function was in $1916$ proposed by an Indian mathematician
Srinivasa Ramanujan. As an analogue of the Riemann hypothesis, an
English mathematician Godfrey Harold Hardy proposed in $1940$ that the real
part of all complex zeros of the Ramanujan zeta function is $6$. This is the
well-known Riemann-Hardy hypothesis for the Ramanujan zeta function.
This article is devoted to the proof of this hypothesis derived from the Ramanujan-Rankin
function.
Owing to the integral representation of the Ramanujan-De
Bruijn function, we establish its series.
We also reduce its product
using the Hadamard's factorization theorem.
By a class with its series and product representations,
 we conclude that the real part of all zeros for Ramanujan-De
Bruijn function is zero.
we also obtain its products of Conrey and Ghosh and Hadamard-type
for the Ramanujan-Rankin function.
Based on the obtained
result, we prove that the Riemann-Hardy hypothesis is true.
\end{abstract}

\maketitle
\tableofcontents

\section{Background and results} {\label{sec:1}}

\subsection{The Ramanujan zeta function}
In 1916 remarkable paper of Srinivasa Ramanujan \cite{1}, who is an Indian
mathematical genius, he proposed the well-known Ramanujan zeta function
$H_\tau \left( x \right)$, defined as
\begin{equation}
\label{eq1}
H_\tau \left( x \right)=\sum\limits_{\ell =1}^\infty {\frac{\tau \left( \ell
\right)}{\ell ^x}} ,
\end{equation}
where $x\in \mathfrak{A}$, $\Re\left( x \right)>13/2$ and $\ell \in \mathfrak{G}$, and
$\tau \left( \ell \right)$ is the Ramanujan's arithmetical function \cite{2}.
Here, let $\mathfrak{A}$, $\mathfrak{B}$ and $\mathfrak{G}$ be the sets of complex, real
and natural numbers, respectively. Suppose $\Re \left( x \right)$ and $\Im
\left( x \right)$ are the real and imaginary parts for a complex variable
$x\in \mathfrak{A}$, respectively.

As the first conjecture of Ramanujan \cite{1}, Mordell \cite{3,4} proved in 1917 that
$\tau \left( \ell \right)$ is multiplicative, i.e.,
\begin{equation}
\label{eq2}
\tau \left( {\ell _1 \ell _2 } \right)=\tau \left( {\ell _1 } \right)\tau
\left( {\ell _2 } \right),
\end{equation}
if $\gcd \left( {\ell _1, \ell _2 } \right)=1$.

Ramanujan \cite{1} said that $\tau \left( \ell \right)$ can also be obtained by
(see \cite{5}, p.161; also see \cite{6,7})
\begin{equation}
\label{eq3}
\begin{array}{l}
L\left( x \right)=q\prod\limits_{n=1}^\infty {\left( {1-q^n} \right)^{24}}
=\sum\limits_{m=1}^\infty {\tau \left( m \right)q^m} ,
\end{array}
\end{equation}
which is the normalized weight 12 cusp form for $SL_2 \left( {\rm Z}
\right)$ \cite{7}, where $L\left( x \right)$ is the modular form, $q=\exp \left(
{2\pi ix} \right)$, $\Im \left( x \right)>0$ and $m,n\in \mathfrak{G}$

As the second conjecture of Ramanujan \cite{1}, Mordell \cite{3} also proved that
\begin{equation}
\label{eq4}
\tau \left( {p^r} \right)=\tau \left( p \right)\tau \left( {p^{r-1}}
\right)-p^{11}\tau \left( {p^{r-2}} \right),
\end{equation}
if $p$ is prime, and $r>2$.

Deligne (see \cite{8}; for the review of the problem, also see \cite{9}) proved that third conjecture of Ramanujan \cite{1}, i.e.,
\begin{equation}
\label{eq5}
\left| {\tau \left( p \right)} \right|\le 2p^{\frac{11}{2}}.
\end{equation}
Ramanujan \cite{1} showed that Eq. (\ref{eq1}) has an Euler-type product of the form \cite{10}
\begin{equation}
\label{eq6}
\begin{array}{l}
H_\tau \left( x \right)=\prod\limits_p {\frac{1}{1-\tau \left( \ell
\right)p^{-x}+p^{11-2x}}} ,
\end{array}
\end{equation}
where $\Re\left( s \right)>13/2$.

\subsection{The Ramanujan-Rankin function}
In 1939 Rankin paper \cite{11}, he first proposed
the well-known Ramanujan-Rankin function $\xi _\tau \left( x \right)$, defined by
\begin{equation}
\label{eq7}
\xi _\tau \left( x \right)=\left( {2\pi } \right)^{-x}\Gamma \left( x
\right)H_\tau \left( x \right),
\end{equation}
where $x\in \mathfrak{A}$ and $\Gamma \left( x
\right)$ is the gamma function.

The modular form $L\left( \omega \right)$ has the functional equation of the form
\cite{11,12}
\begin{equation}
\label{eq8}
L\left( \omega \right)=\frac{1}{\omega ^{12}}L\left( {-\frac{1}{\omega }}
\right)
\end{equation}
for $\omega >0$.

Let $i=\sqrt {-1} $. There exists an integral of (\ref{eq7}) as follows \cite{11,12}:
\begin{equation}
\label{eq9}
\xi _\tau \left( x \right)=\int\limits_0^\infty {\omega ^{x-1}L\left(
{i\omega } \right)d\omega } ,
\end{equation}
which can be rewritten as (\cite{11})
\begin{equation}
\label{eq10}
\xi _\tau \left( x \right)=\int\limits_0^1 {L\left( {i\omega } \right)\left(
{\omega ^{x-1}+\omega ^{11-x}} \right)d\omega } .
\end{equation}
Wilton \cite{13} suggested an alternative integral representation of (\ref{eq7}), given as \cite{13a}
\begin{equation}
\label{eq11}
\xi _\tau \left( x \right)=\int\limits_0^\infty {\omega ^{x-1}{\rm L}\left(
{\exp \left( {-2\pi \omega } \right)} \right)d\omega } =\int\limits_1^\infty
{\left( {\omega ^{x-1}+\omega ^{11-x}} \right){\rm L}\left( {\exp \left(
{-2\pi \omega } \right)} \right)d\omega } ,
\end{equation}
provided that
\begin{equation}
\label{eq12}
\begin{array}{l}
{\rm L}\left( q \right)=q\prod\limits_{n=1}^\infty {\left( {1-q^n}
\right)^{24}} =\sum\limits_{m=1}^\infty {\tau \left( m \right)q^m} .
\end{array}
\end{equation}
Rankin \cite{11} and Wilton \cite{13} proved that (\ref{eq7}) has the functional equation
as follows (also see \cite{12}):
\begin{equation}
\label{eq13}
\xi _\tau \left( x \right)=\xi _\tau \left( {12-x} \right).
\end{equation}

\subsection{The product of Conrey and Ghosh}

Conrey and Ghosh \cite{14} suggested the product of (\ref{eq7}) as follows:
\begin{equation}
\label{eq14}
\begin{array}{l}
\xi _\tau \left( x \right)=e^{b+b_0 x}\prod\limits_{x_m } {\left(
{1-\frac{x}{x_m }} \right)\exp \left( {\frac{x}{x_m }} \right)} ,
\end{array}
\end{equation}
where $x_m $ take over all zeros of (\ref{eq7}),
and both $b$ and $b_0 $ are two unknown constants.
Here, (\ref{eq14}) is so-called product of Conrey and Ghosh.

Wilton \cite{13} considered $x=6+i\varphi $ into (\ref{eq13}) and obtained the
Ramanujan-Wilton function $\Lambda \left( \varphi \right)$, defined by
\begin{equation}
\label{eq15}
\Lambda \left( \varphi \right)=\xi _\tau \left( {6+i\varphi } \right)=\xi
_\tau \left( {6-i\varphi } \right),
\end{equation}
which satisfies the functional equation \cite{13,15}
\begin{equation}
\label{eq16}
\Lambda \left( {-\varphi } \right)=\Lambda \left( {\varphi } \right),
\end{equation}
where $\varphi\in \mathfrak{A}$.

In 1950 de Bruijn \cite{15} presented the Ramanujan-de Bruijn function $X\left( z \right)$, defined as
\begin{equation}
\label{eq17}
X\left( z \right):=\Lambda \left( {iz} \right)=\left( {2\pi }.
\right)^{-\left( {6+z} \right)}\Gamma \left( {6+z} \right)H_\tau \left(
{6+z} \right),
\end{equation}
where $z\in \mathfrak{A}$.

\subsection{The Riemann-Hardy hypothesis}
In analogy with the Riemann hypothesis, Godfrey Harold Hardy (see \cite{15}, p.174),
who is an English mathematician, proposed the well-known conjecture that
states the following:

\begin{theorem}
\label{TH1}
(Riemann-Hardy hypothesis)

All of the nontrivial zeros of the Ramanujan zeta function (\ref{eq1}) lie on the
critical line $\Re \left( x \right)=6$.
\end{theorem}

Meanwhile, Hardy (see \cite{15}, p.174) also conjectured an
analogue of the Rieman-Siegel-like formula and von Mangoldt-like theorem. As
an analogous work of Hardy theorem for the Riemann zeta function, it was
proved by Lekkerkerker \cite{18}. The latter was studied by generalized by Berndt
\cite{19} and further developed by Hafner \cite{20}, Ki \cite{21} and Chirre and
Casta\~{n}\'{o}n \cite{22}. The former was developed by Ferguson and coauthors
\cite{23} and further reported by Keiper \cite{24}. More importantly, the papers of
Ferguson and coauthors \cite{23} and Keiper \cite{24} reported some nontrivial zeros for (\ref{eq1}).
As one of interesting and important problems in the theory of modular form
\cite{25}, the Riemann-Hardy hypothesis has not been solved
for the more than nineteen-eighty years.

Wilton \cite{13} reported (\ref{eq7}) is an entire function,
and
Apostol and Sklar \cite{26} said that (\ref{eq7}) is an entire
function of finite order.
Thanks to
the result of Ogg (see \cite{26a}, p.21), it is proved that (\ref{eq7}) is
an entire function of order $\lambda =1$.
There is an unsolved problem that the product of Conrey and Ghosh
\cite{14} has the unknown parameters.
The series and  product of Hadamard-type representations of
(\ref{eq7}) has not considered in any paper.
Theorem \ref{TH1} implies that
the real
part of all complex zeros of (\ref{eq7}) is $6$.

In fact, de Bruijn \cite{15} guessed that
Theorem \ref{TH1} is equivalent to the case that
all zeros of (\ref{eq17}) lie on the critical line $\Re \left(
z \right)=0$. This is the conjecture of de Bruijn.
Moreover, there is no investigation for the series, product, zeros and order of
(\ref{eq17}) at present. The conjecture of de Bruijn remains an open problem.

There exists an equivalent representation of
the reported works, which implies that
all zeros of (\ref{eq15}) lie on the critical
line $\Im \left( \varphi \right)=0$.
This is the conjecture of Ki \cite{20},
which is further studied by Chirre and
Casta\~{n}\'{o}n \cite{21}.
The series, product and order of (\ref{eq15})
have not reported. The conjecture of Ki is also an open problem.

\subsection{The main aim of the present paper}

Motivated by the ideas above, we
give the outline of the present paper. In Section 2 we present the order,
the series and product presentations of (\ref{eq7}), (\ref{eq15}) and (\ref{eq17}).
The values of the parameters in the product of Conrey and Ghosh are given in detail.
The Hadamard-type product are also obtained. The product for (\ref{eq17})
is set up based on the functional equation
and Hadamard¨s factorization theorem.
In
Section 3 we first consider the detailed proof of the conjecture of de Bruijn.
We then reduce to the conjecture of Ki.
Finally, we report the proof of Theorem \ref{TH1}.

\section{Preliminary results}

\subsection{The series and orders }
To begin with, we investigate the series representation of $\xi _\tau \left(
x \right)$ as follows.

Let $\mathfrak{H}=\mathfrak{G}\cup \left\{ 0 \right\}$.

\begin{proposition}
\label{PRO1}
Let $x\in \mathfrak{A}$. Then there exists
\begin{equation}
\label{eq18}
\xi _\tau \left( x \right)=\sum\limits_{s=0}^\infty {\chi \left( s
\right)\left( {x-6} \right)^{2s}}
\end{equation}
where
\begin{equation}
\label{eq19}
\chi \left( s \right)=2\int\limits_0^1 {L\left( {i\omega } \right)\omega
^5\left[ {\frac{\left( {\ln \omega } \right)^{2s}}{\left( {2s} \right)!}}
\right]d\omega }
\end{equation}
are the positive real coefficients and $s\in \mathfrak{H}$.
\end{proposition}

\begin{proof}
By using (\ref{eq11}), we present
\begin{equation}
\label{eq20}
 \begin{aligned}
 \xi _\tau \left( x \right)&=\int\limits_0^1 {L\left( {i\omega }
\right)\left( {\omega ^{x-1}+\omega ^{11-x}} \right)d\omega } \\
&=\int\limits_0^1 {L\left( {i\omega } \right)\omega ^5\left( {\omega
^{x-6}+\omega ^{6-x}} \right)d\omega } \\
&=\int\limits_0^1 {L\left( {i\omega } \right)\omega ^5\left\{ {\exp \left[
{\left( {x-6} \right)\ln \omega } \right]+\exp \left[ {\left( {6-x}
\right)\ln \omega } \right]} \right\}d\omega } . \\
 \end{aligned}
\end{equation}
To simply (\ref{eq20}), we obtain
\begin{equation}
\label{eq21}
\xi _\tau \left( x \right)=2\int\limits_0^1 {L\left( {i\omega }
\right)\omega ^5\cosh \left[ {\left( {x-6} \right)\ln \omega }
\right]d\omega } .
\end{equation}
Consider that
\begin{equation}
\label{eq22}
\cosh \left[ {\left( {x-6} \right)\ln \omega }
\right]=\sum\limits_{s=0}^\infty {\frac{\left[ {\left( {x-6} \right)\ln
\omega } \right]^{2s}}{\left( {2s} \right)!}} .
\end{equation}
By the substitution of (\ref{eq22}) into (\ref{eq21}),
\begin{equation}
\label{eq23}
 \begin{aligned}
 \xi _\tau \left( x \right)&=2\int\limits_0^1 {L\left( {i\omega }
\right)\omega ^5\left\{ {\cosh \left[ {\left( {x-6} \right)\ln \omega }
\right]} \right\}d\omega } \\
&=2\int\limits_0^1 {L\left( {i\omega } \right)\omega ^5\left\{
{\sum\limits_{s=0}^\infty {\frac{\left[ {\left( {x-6} \right)\ln \omega }
\right]^{2s}}{\left( {2s} \right)!}} } \right\}d\omega }\\
&=\sum\limits_{s=0}^\infty {\left\{ {2\int\limits_0^1 {L\left( {i\omega }
\right)\omega ^5\left[ {\frac{\left( {\ln \omega } \right)^{2s}}{\left( {2s}
\right)!}} \right]d\omega } } \right\}\left( {x-6} \right)^{2s}} . \\
 \end{aligned}
\end{equation}
In fact,
\begin{equation}
\label{eq24}
\chi \left( s \right)=2\int\limits_0^1 {L\left( {i\omega } \right)\omega
^5\left[ {\frac{\left( {\ln \omega } \right)^{2s}}{\left( {2s} \right)!}}
\right]d\omega } >0
\end{equation}
provided that
\begin{equation}
\label{eq25}
\begin{array}{l}
L\left( {i\omega } \right)=\exp \left( {-2\pi \omega }
\right)\prod\limits_{n=1}^\infty {\left[ {1-\exp \left( {-2n\pi \omega }
\right)} \right]^{24}} >0
\end{array}
\end{equation}
is true.

We thus complete the proof.
\end{proof}

\begin{proposition}
\label{PRO2}
If $\chi \left( s \right)$ is defined as in Proposition \ref{PRO1}, then
there exists
\begin{equation}
\label{eq26}
\Lambda \left( \varphi \right)=\sum\limits_{s=0}^\infty {\chi \left( s
\right)\left( {-1} \right)^s\varphi ^{2s}} ,
\end{equation}
where $\varphi \in \mathfrak{A}$.
\end{proposition}

\begin{proof}
Putting $x=6+i\varphi $ in Proposition \ref{PRO1} implies the desired
result.
\end{proof}

 Here, we need to prove the following:

\begin{proposition}
\label{PRO3}
$\Lambda \left( \varphi \right)$
is an even entire function of order $\lambda=1$.
\end{proposition}

\begin{proof}
Owing to (\ref{eq26}), we rewrite
\begin{equation}
\label{eq27}
\chi \left( s \right)=2\int\limits_0^1 {L\left( {i\omega } \right)\omega
^5\left[ {\frac{\left( {\ln \omega } \right)^{2s}}{\left( {2s} \right)!}}
\right]d\omega }
\end{equation}
as
\begin{equation}
\label{eq28}
\chi \left( s \right)=\frac{2}{\left( {2s} \right)!}\int\limits_0^1 {L\left(
{i\omega } \right)\omega ^5\exp \left( {2s\ln \omega } \right)d\omega } >0.
\end{equation}
By virtue of
\begin{equation}
\label{eq29}
0< \exp \left( {2s\log \omega } \right)\le 1
\end{equation}
for $0< \omega \le 1$, there exists
\begin{equation}
\label{eq30}
\left| {\chi \left( s \right)} \right|=\frac{2}{\left( {2s} \right)!}
{\int\limits_0^1 {L\left( {i\omega } \right)\omega ^5\exp \left( {2s\ln
\omega } \right)d\omega } } .
\end{equation}
Because of
\begin{equation}
\label{eq31}
0< {\int\limits_0^1 {L\left( {i\omega } \right)\omega ^5\exp \left(
{2s\ln \omega } \right)d\omega } } \le C,
\end{equation}
where
\begin{equation}
\label{eq32}
C= {\int\limits_0^1 {L\left( {i\omega } \right)\omega ^5d\omega } }
=\int\limits_0^1 {L\left( {i\omega } \right)\omega ^5d\omega } =\xi
_\tau \left( 6 \right)<\infty ,
\end{equation}
we have from (\ref{eq30}) and (\ref{eq31}) that
\begin{equation}
\label{eq33}
\left| {\chi \left( s \right)} \right|=\frac{2}{\left( {2s} \right)!}
{\int\limits_0^1 {L\left( {i\omega } \right)\omega ^5\exp \left( {2s\ln
\omega } \right)d\omega } } \le \frac{2C}{\left( {2s} \right)!}.
\end{equation}
By the formula (1.04) in Levin' book (\cite{27}, p.4), we have
\begin{equation}
\label{eq34}
\mathop {\lim }\limits_{s\to \infty } \sqrt[s]{\left| {\chi \left( s
\right)} \right|}=0
\end{equation}
such that $\Lambda \left( \varphi \right)$ is an entire function.

Using Ogg's result (see \cite{26a}, p.21) that
$\xi _\tau \left( x \right)$ is an entire function of order $\lambda=1$,
Theorem 2 in Levin' book (\cite{27}, p.4) said that the order $\lambda $ of
$\Lambda \left( \varphi \right)$ can be also given as
\begin{equation}
\label{eq39}
\lambda =\overline {\mathop {\lim }\limits_{s\to \infty } } \frac{s\ln
s}{\ln \frac{1}{\left| {\chi \left( s \right)} \right|}}=1.
\end{equation}
By (\ref{eq39}), $\Lambda \left(
\varphi \right)$ is also an entire function of
order $\lambda =1$ because $|\chi \left( s
\right)\left( {-1} \right)^s|=\chi \left( s
\right)>0$.

We thus finish the proof.
\end{proof}

Putting $x=6+z$ in Proposition \ref{PRO1} for $z\in \mathfrak{A}$ and using the
formulae (\ref{eq9}) and (\ref{eq11}) of Rankin \cite{11}, we obtain
\begin{equation}
\label{eq42}
 \begin{aligned}
 X\left( z \right)&=\left( {2\pi } \right)^{-\left( {6+z} \right)}\Gamma
\left( {6+z} \right)H_\tau \left( {6+z} \right)=\int\limits_0^\infty {\omega
^{5+z}L\left( {i\omega } \right)d\omega } \\
&=\int\limits_0^1 {L\left( {i\omega } \right)\omega ^5\left( {\omega
^z+\omega ^{-z}} \right)d\omega }, \\
 \end{aligned}
\end{equation}
which is the result of De Bruijn \cite{15}.

Based on this, we get the followings:

\begin{corollary}
\label{CORO1}
If $\chi \left( s \right)$ is defined as in Proposition \ref{PRO1}, then we
have
\begin{equation}
\label{eq43}
X\left( z \right)=\sum\limits_{s=0}^\infty {\chi \left( s \right)z^{2s}} ,
\end{equation}
where $z\in \mathfrak{A}$.
\end{corollary}

\begin{proof}
As a similar manner of Proposition \ref{PRO1}, we have from (\ref{eq42}) the
required result.
\end{proof}

Moreover, we have the followings:

\begin{corollary}
\label{CORO2}
The function $X\left( z \right)$ is an even entire function of order
$\lambda =1$.
\end{corollary}

\begin{proof}
With (\ref{eq43}), we know
\begin{equation}
\label{eq44}
X\left( z \right)=X\left( {-z} \right),
\end{equation}
where $z \in \mathfrak{A}$, and by using (\ref{eq34}) and (\ref{eq39}), we see that $X\left( z \right)$ is an entire
function of order $\lambda =1$.

Hence, we complete the proof.
\end{proof}

\subsection{The products}

\begin{theorem}
\label{TH3}
There exists
\begin{equation}
\label{eq46}
\begin{array}{l}
X\left( z \right)=X\left( 0 \right)\prod\limits_{\Im \left( {z_m }
\right)>0} {\left( {1-\frac{z^2}{z_m^2 }} \right)} ,
\end{array}
\end{equation}
where the product run over all of the zeros of $z_m $ of $X\left( z
\right)$.
\end{theorem}

\begin{proof}

Taking $z=0$ in (\ref{eq42}) implies
\begin{equation}
\label{eq47}
X\left( 0 \right)=\int\limits_0^\infty {\omega ^5L\left( {i\omega }
\right)d\omega } >0
\end{equation}
provided that (\ref{eq25}) holds.

In view of Corollary \ref{CORO2}, the Hadamard's factorization theorem (see
\cite{28}, p.250) said that
\begin{equation}
\label{eq48}
\begin{array}{l}
X\left( z \right)=e^{M+Lz}\prod\limits_{z_m } {\left( {1-\frac{z}{z_m }}
\right)\exp \left( {\frac{z}{z_m }} \right)} ,
\end{array}
\end{equation}
where the product run over all of the zeros $z_m $ of $X\left( z \right)$
and both $M$ and $L$ are two constants.
Clearly, the sequence $\left\{ {z_m } \right\}$ are the set of non-zero zeros of $X\left( z \right)$ with $|z _{m+1}|>|z _{m}|$.

Since $X\left( z \right)$ is an even entire function obtained in Corollary
\ref{CORO2},
\begin{equation}
\label{eq49}
 \begin{aligned}
 X\left( z \right)&=e^{M+Lz}\prod\limits_{z_m } {\left( {1-\frac{z}{z_m }}
\right)\exp \left( {\frac{z}{z_m }} \right)} \\
&=e^{M+Lz}\prod\limits_{\Im \left( {z_m } \right)>0} {\left( {1-\frac{z}{z_m
}} \right)\left( {1+\frac{z}{z_m }} \right)\exp \left( {\frac{z}{z_m
}-\frac{z}{z_m }} \right)} \\
&=e^{M+Lz}\prod\limits_{\Im \left( {z_m }
\right)>0} {\left( {1-\frac{z^2}{z_m^2 }} \right)} , \\
 \end{aligned}
\end{equation}
which yields that
\begin{equation}
\label{eq50}
\begin{array}{l}
X\left( {-z} \right)=e^{M-Lz}\prod\limits_{\Im \left( {z_m } \right)>0}
{\left( {1-\frac{z^2}{z_m^2 }} \right)} .
\end{array}
\end{equation}
Combining (\ref{eq49}) and (\ref{eq50}) and using (\ref{eq44}), we suggest
\begin{equation}
\label{eq51}
\begin{array}{l}
e^{M+Lz}\prod\limits_{\Im \left( {z_m } \right)>0} {\left(
{1-\frac{z^2}{z_m^2 }} \right)} =e^{M-Lz}\prod\limits_{\Im \left( {z_m }
\right)>0} {\left( {1-\frac{z^2}{z_m^2 }} \right)} .
\end{array}
\end{equation}
By (\ref{eq47}) and (\ref{eq51}), we have $M=\ln X\left( 0 \right)$ and $L=0$ such that
(\ref{eq49}) can be rewritten as
\begin{equation}
\label{eq52}
\begin{array}{l}
X\left( z \right)=X\left( 0 \right)\prod\limits_{\Im \left( {z_m }
\right)>0} {\left( {1-\frac{z^2}{z_m^2 }} \right)} .
\end{array}
\end{equation}
Thus, we obtain the required result.
\end{proof}

As a direct result, we obtain:

\begin{corollary}
\label{CORO4}
There exists any positive number $\varepsilon >0$ such that
\begin{equation}
\label{eq53}
\frac{1}{\left| {z_m } \right|^{1+\varepsilon }}
\end{equation}
is convergent.
\end{corollary}

\begin{proof}
Since $X\left( z \right)$ is an even entire function of order $\lambda =1$
and (\ref{eq46}) is true, as an analogous result in Levin' book (\cite{27}, p.8), we
have (\ref{eq53}) such that $\sum\limits_{m=1}^\infty {\left| {z/z_m }
\right|^{1+\varepsilon }} $ converges uniformly in each bounded domain.

This is the required result.
\end{proof}

\subsubsection{The product of Conrey and Ghosh}

As a direct result of Theorem \ref{TH3}, we see the following:

\begin{theorem}
\label{TH4}
Assume the above denotations. Then we have:
\begin{itemize}
  \item
(B1) There is
\begin{equation}
\label{eq54}
\begin{array}{l}
\xi _\tau \left( x \right)=\xi _\tau \left( 6 \right)\prod\limits_{\Im
\left( {z_m } \right)>0} {\left[ {1-\frac{\left( {x-6} \right)^2}{z_m^2 }}
\right]} .
\end{array}
\end{equation}
  \item
(B2) (\textbf{The product of Conrey and Ghosh}) There is
\begin{equation}
\label{eq55}
\begin{array}{l}
\xi _\tau \left( x \right)=e^{b+b_0 x}\prod\limits_{x_m } {\left(
{1-\frac{x}{x_m }} \right)\exp \left( {\frac{x}{x_m }} \right)} ,
\end{array}
\end{equation}
where the product takes over all of the zeros $x_m $ of $\xi _\tau \left( x
\right)$, $b_0 =\ln \xi _\tau \left( 0 \right)-6\sum\limits_{x_m } {\left(
{1/x_m } \right)} $ and $b=\ln \xi _\tau \left( 0 \right)$.
  \item
(B3) (\textbf{The product of Hadamard-type}) There is
\begin{equation}
\label{eq56}
\begin{array}{l}
\xi _\tau \left( x \right)=\xi _\tau \left( 0 \right)\prod\limits_{x_m }
{\left( {1-\frac{x}{x_m }} \right)} .
\end{array}
\end{equation}
Moreover, there exists any positive number $\delta >0$ such that
$\sum\limits_{m=1}^\infty {\left| {x_m } \right|^{-\left( {1+\delta }
\right)}} $ is convergent.
\end{itemize}
\end{theorem}

\begin{proof}
Step one is to prove (B1). Taking $z=x-6$ in Theorem \ref{TH3} and Corollary \ref{CORO1}, we
have
\begin{equation}
\label{eq57}
\begin{array}{l}
X\left( {x-6} \right)=X\left( 0 \right)\prod\limits_{\Im \left( {z_m }
\right)>0} {\left[ {1-\frac{\left( {x-6} \right)^2}{z_m^2 }} \right]}
\end{array}
\end{equation}
and
\begin{equation}
\label{eq58}
X\left( {x-6} \right)=\sum\limits_{s=0}^\infty {\chi \left( s \right)\left(
{x-6} \right)^{2s}} .
\end{equation}
By Proposition \ref{PRO1}, we see that (\ref{eq58}) is equal to
\begin{equation}
\label{eq59}
\begin{array}{l}
\xi _\tau \left( x \right)=X\left( {x-6} \right)=\sum\limits_{s=0}^\infty
{\chi \left( s \right)\left( {x-6} \right)^{2s}} =X\left( 0
\right)\prod\limits_{\Im \left( {z_m } \right)>0} {\left[ {1-\frac{\left(
{x-6} \right)^2}{z_m^2 }} \right]} .
\end{array}
\end{equation}
Because (\ref{eq59}) gives
\[
\xi _\tau \left( 6 \right)=X\left( 0 \right),
\]
(\ref{eq59}) becomes
\begin{equation}
\label{eq60}
\begin{array}{l}
\xi _\tau \left( x \right)=\xi _\tau \left( 6 \right)\prod\limits_{\Im
\left( {z_m } \right)>0} {\left[ {1-\frac{\left( {x-6} \right)^2}{z_m^2 }}
\right]} .
\end{array}
\end{equation}
Thus, we finish the proof of (B1).

Step two is to prove (B2).

By Proposition \ref{PRO3}, we adopt the Hadamard's factorization theorem in
Titchmarsh's book (\cite{28}, p.250) to obtain the conjecture of Conrey and Ghosh
\cite{14}, that is,
\begin{equation}
\label{eq61}
\begin{array}{l}
\xi _\tau \left( x \right)=e^{b+b_0 x}\prod\limits_{x_m } {\left(
{1-\frac{x}{x_m }} \right)\exp \left( {\frac{x}{x_m }} \right)} ,
\end{array}
\end{equation}
where the product run over all of the zeros $x_m $ of $\xi _\tau \left( x \right)$
and both $b$ and $b_0 $ are two constants.
Obviously, the sequence $\left\{ {x_m } \right\}$ is the set of non-zero zeros of $\xi _\tau \left( x \right)$ with
$|x _{m+1}|>|x _{m}|$.

Combining (\ref{eq21}), (\ref{eq60}) and (\ref{eq61}), we have
\begin{equation}
\label{eq62}
 \begin{aligned}
 \xi _\tau \left( x \right)&=e^{b+b_0 x}\prod\limits_{x_m } {\left(
{1-\frac{x}{x_m }} \right)\exp \left( {\frac{x}{x_m }} \right)} \\
&=\xi _\tau \left( 6 \right)\prod\limits_{\Im \left( {z_m } \right)>0}
{\left[ {1-\frac{\left( {x-6} \right)^2}{z_m^2 }} \right]} \\
&=2\int\limits_0^1 {L\left( {i\omega } \right)\omega ^5\cosh \left[ {\left(
{x-6} \right)\ln \omega } \right]d\omega } \\
 \end{aligned}
\end{equation}
such that
\begin{equation}
\label{eq63}
\left( {x_m -6} \right)^2=z_m^2 .
\end{equation}
From (\ref{eq62}) and (\ref{eq63}) we may get
\begin{equation}
\label{eq64}
\begin{array}{l}
\xi _\tau \left( x \right)=\xi _\tau \left( 6 \right)\prod\limits_{\Im
\left( {z_m } \right)>0} {\left[ {1-\frac{\left( {x-6} \right)^2}{\left(
{x_m -6} \right)^2}} \right]} .
\end{array}
\end{equation}
In view of the fact $x_m =6\pm z_m $ obtained by (\ref{eq63}), we have
\begin{equation}
\label{eq65}
\Im \left( {z_m } \right)=\Im \left( {x_m -6} \right)=\Im \left( {x_m }
\right)
\end{equation}
such that (\ref{eq62}) becomes
\begin{equation}
\label{eq66}
 \begin{aligned}
 \xi _\tau \left( x \right)&=\xi _\tau \left( 6 \right)\prod\limits_{\Im
\left( {z_m } \right)>0} {\left[ {1-\frac{\left( {x-6} \right)^2}{\left(
{x_m -6} \right)^2}} \right]} \\
&=\xi _\tau \left( 6 \right)\prod\limits_{\Im
\left( {x_m -6} \right)>0} {\left[ {1-\frac{\left( {x-6} \right)^2}{\left(
{x_m -6} \right)^2}} \right]} \\
 &=\xi _\tau \left( 6 \right)\prod\limits_{\Im \left( {x_m } \right)>0}
{\left[ {1-\frac{\left( {x-6} \right)^2}{\left( {x_m -6} \right)^2}}
\right]} . \\
 \end{aligned}
\end{equation}
Let us write
\begin{equation}
\label{eq67}
 \begin{aligned}
 \xi _\tau \left( x \right)&=\xi _\tau \left( 6 \right)\prod\limits_{\Im
\left( {x_m } \right)>0} {\left[ {1-\left( {\frac{x-6}{x_m -6}} \right)^2}
\right]} \\
&=\xi _\tau \left( 6 \right)\prod\limits_{\Im \left( {x_m } \right)>0}
{\left\{ {\left( {1-\frac{x-6}{x_m -6}} \right)\left( {1+\frac{x-6}{x_m -6}}
\right)} \right\}} \\
&=\xi _\tau \left( 6 \right)\prod\limits_{\Im \left( {x_m } \right)>0}
{\left\{ {\left( {1-\frac{x-6}{x_m -6}} \right)\left( {1-\frac{x-6}{6-x_m }}
\right)} \right\}} \\
&=\xi _\tau \left( 6 \right)\prod\limits_{\Im \left( {x_m } \right)>0}
{\left\{ {\left( {1-\frac{x-6}{x_m -6}} \right)\left( {1-\frac{x-6}{6-x_m }}
\right)} \right\}} \\
&=\xi _\tau \left( 6 \right)\prod\limits_{\Im \left( {x_m } \right)>0}
{\left( {1-\frac{x-6}{x_m -6}} \right)\left[ {1-\frac{x-6}{\left( {12-x_m }
\right)-6}} \right]} . \\
 \end{aligned}
\end{equation}
By using the functional equation (\ref{eq13}), the identity (\ref{eq67}) is equal to
\begin{equation}
\label{eq68}
 \begin{aligned}
 \xi _\tau \left( x \right)&=\xi _\tau \left( 6 \right)\prod\limits_{\Im
\left( {x_m } \right)>0} {\left( {1-\frac{x-6}{x_m -6}} \right)\left[
{1-\frac{x-6}{\left( {12-x_m } \right)-6}} \right]} \\
&=\xi _\tau \left( 6 \right)\prod\limits_{x_m } {\left( {1-\frac{x-6}{x_m
-6}} \right)} . \\
 \end{aligned}
\end{equation}
From (\ref{eq62}) we take into account
\begin{equation}
\label{eq69}
\begin{array}{l}
\xi _\tau \left( x \right)=e^{b+b_0 x}\prod\limits_{x_m } {\left(
{1-\frac{x}{x_m }} \right)\exp \left( {\frac{x}{x_m }} \right)}
=2\int\limits_0^1 {L\left( {i\omega } \right)\omega ^5\cosh \left[ {\left(
{x-6} \right)\ln \omega } \right]d\omega }
\end{array}
\end{equation}
such that
\begin{equation}
\label{eq70}
\xi _\tau \left( 0 \right)=e^b=2\int\limits_0^1 {L\left( {i\omega }
\right)\omega ^5\cosh \left( {6\ln \omega } \right)d\omega }
=2\int\limits_0^1 {L\left( {i\omega } \right)\omega ^{11}d\omega } >0,
\end{equation}
\begin{equation}
\label{eq71}
\begin{array}{l}
\xi _\tau \left( 6 \right)=e^{b+6b_0 }\prod\limits_{x_m } {\left(
{1-\frac{6}{x_m }} \right)\exp \left( {\frac{6}{x_m }} \right)}
=2\int\limits_0^1 {L\left( {i\omega } \right)\omega ^5d\omega } >0
\end{array}
\end{equation}
and
\begin{equation}
\label{eq72}
 \begin{aligned}
 \xi _\tau \left( {12} \right)&=e^{b+12b_0 }\prod\limits_{x_m } {\left(
{1-\frac{12}{x_m }} \right)\exp \left( {\frac{12}{x_m }} \right)}\\
&=2\int\limits_0^1 {L\left( {i\omega } \right)\omega ^5\cosh \left( {6\ln
\omega } \right)d\omega } \\
&=2\int\limits_0^1 {L\left( {i\omega } \right)\omega ^{11}d\omega } >0. \\
 \end{aligned}
\end{equation}
Making use of (\ref{eq70}), (\ref{eq71}) and (\ref{eq72}), we have
\begin{equation}
\label{eq73}
b=\ln \xi _\tau \left( 0 \right)
\end{equation}
and
\begin{equation}
\label{eq74}
\xi _\tau \left( 0 \right)=\xi _\tau \left( {12} \right)
\end{equation}
such that
\begin{equation}
\label{eq75}
\begin{array}{l}
\xi _\tau \left( 0 \right)=e^b=e^{b+12b_0 }\prod\limits_{x_m } {\left(
{1-\frac{12}{x_m }} \right)\exp \left( {\frac{12}{x_m }} \right)}
=2\int\limits_0^1 {L\left( {i\omega } \right)\omega ^{11}d\omega } >0.
\end{array}
\end{equation}
By using (\ref{eq68}), we obtain
\begin{equation}
\label{eq76}
\begin{array}{l}
\xi _\tau \left( x \right)=\xi _\tau \left( 6 \right)\prod\limits_{x_m }
{\left( {1-\frac{x-6}{x_m -6}} \right)} ,
\end{array}
\end{equation}
which leads to
\begin{equation}
\label{eq77}
\begin{array}{l}
\xi _\tau \left( 0 \right)=\xi _\tau \left( 6 \right)\prod\limits_{x_m }
{\left( {1-\frac{0-6}{x_m -6}} \right)} =\xi _\tau \left( 6
\right)\prod\limits_{x_m } {\left( {1+\frac{6}{x_m -6}} \right)} >0.
\end{array}
\end{equation}
From (\ref{eq77}) we have
\begin{equation}
\label{eq78}
\begin{array}{l}
\xi _\tau \left( 0 \right)=\xi _\tau \left( 6 \right)\prod\limits_{x_m }
{\left( {\frac{x_m }{x_m -6}} \right)}
\end{array}
\end{equation}
such that
\begin{equation}
\label{eq79}
\begin{array}{l}
\xi _\tau \left( 6 \right)=\xi _\tau \left( 0 \right)\prod\limits_{x_m }
{\left( {\frac{x_m -6}{x_m }} \right)} =\xi _\tau \left( 0
\right)\prod\limits_{x_m } {\left( {1-\frac{6}{x_m }} \right)} .
\end{array}
\end{equation}
From (\ref{eq71}) and (\ref{eq73}) we obtain
\begin{equation}
\label{eq80}
\begin{array}{l}
\xi _\tau \left( 6 \right)=\xi _\tau \left( 0 \right)e^{6b_0
}\prod\limits_{x_m } {\left( {1-\frac{6}{x_m }} \right)\exp \left(
{\frac{6}{x_m }} \right)} .
\end{array}
\end{equation}
Combining (\ref{eq79}) and (\ref{eq80}) gives
\begin{equation}
\label{eq81}
\begin{array}{l}
\xi _\tau \left( 6 \right)=\xi _\tau \left( 0 \right)\prod\limits_{x_m }
{\left( {1-\frac{6}{x_m }} \right)} =\left[ {e^{6b_0 }\prod\limits_{x_m }
{\exp \left( {\frac{6}{x_m }} \right)} } \right]\prod\limits_{x_m } {\left(
{1-\frac{6}{x_m }} \right)} ,
\end{array}
\end{equation}
which leads to
\begin{equation}
\label{eq82}
\begin{array}{l}
\xi _\tau \left( 0 \right)=e^{6b_0 }\prod\limits_{x_m } {\exp \left(
{\frac{6}{x_m }} \right)} .
\end{array}
\end{equation}
It follows from (\ref{eq82}) that
\begin{equation}
\label{eq83}
\begin{array}{l}
e^{6b_0 }=\xi _\tau \left( 0 \right)\prod\limits_{x_m } {\exp \left(
{-\frac{6}{x_m }} \right)} ,
\end{array}
\end{equation}
which yields that
\begin{equation}
\label{eq84}
b_0 =\ln \xi _\tau \left( 0 \right)-6\sum\limits_{x_m } {\frac{1}{x_m }} .
\end{equation}
Thus, (\ref{eq61}) is true under the conditions (\ref{eq73}) and (\ref{eq84}).

Step three is to prove (B3).

By using (\ref{eq79}), we reconsider
\begin{equation}
\label{eq85}
\begin{array}{l}
\xi _\tau \left( x \right)=\xi _\tau \left( 6 \right)\prod\limits_{x_m }
{\left( {1-\frac{x-6}{x_m -6}} \right)} =\xi _\tau \left( 0
\right)\prod\limits_{x_m } {\left( {1-\frac{6}{x_m }} \right)}
\prod\limits_{x_m } {\left( {1-\frac{x-6}{x_m -6}} \right)} .
\end{array}
\end{equation}
Further, (\ref{eq85}) becomes
\begin{equation}
\label{eq86}
 \begin{aligned}
 \xi _\tau \left( x \right)&=\xi _\tau \left( 0 \right)\prod\limits_{x_m }
{\left( {1-\frac{6}{x_m }} \right)} \prod\limits_{x_m } {\left(
{1-\frac{x-6}{x_m -6}} \right)} \\
 &=\xi _\tau \left( 0 \right)\prod\limits_{x_m } {\left( {1-\frac{6}{x_m }}
\right)} \cdot \prod\limits_{x_m } {\frac{x_m -6-\left( {x-6} \right)}{x_m
-6}} \\
&=\xi _\tau \left( 0 \right)\prod\limits_{x_m } {\frac{x_m -6}{x_m }} \cdot
\prod\limits_{x_m } {\frac{x_m -s}{x_m -6}} \\
&=\xi _\tau \left( 0
\right)\prod\limits_{x_m } {\left( {\frac{x_m -6}{\rho _k }\cdot \frac{x_m
-s}{x_m -6}} \right)} \\
&=\xi _\tau \left( 0 \right)\prod\limits_{x_m } {\left( {\frac{x_m -6}{x_m
-6}\cdot \frac{x_m -s}{x_m }} \right)} \\
&=\xi _\tau \left( 0
\right)\prod\limits_{x_m } {\frac{x_m -x}{x_m }} . \\
 \end{aligned}
\end{equation}
It follows from (\ref{eq86}) that
\begin{equation}
\label{eq87}
\begin{array}{l}
\xi _\tau \left( x \right)=\xi _\tau \left( 0 \right)\prod\limits_{x_m }
{\frac{x_m -x}{x_m }} =\xi _\tau \left( 0 \right)\prod\limits_{x_m } {\left(
{1-\frac{x}{x_m }} \right)} .
\end{array}
\end{equation}
Clearly, (\ref{eq87}) is the same as (\ref{eq56}).
Thus, this proof is finished.
\end{proof}

By using Proposition \ref{PRO3}, the order $\lambda =1$ of $\xi _\tau \left(
x \right)$ implies that there is any positive number $\delta >0$ such that
$\left| {x_m } \right|^{-\left( {1+\delta } \right)}<\infty $. Thus,
the series $\sum\limits_{m=1}^\infty {\left| {x_m } \right|^{-\left( {1+\delta }
\right)}} $ is convergent.

\begin{remark}
With (\ref{eq21}) and (\ref{eq56}), we establish the relation
\begin{equation}
\label{eq88}
\begin{array}{l}
\xi _\tau \left( x \right)=2\int\limits_0^1 {L\left( {i\omega }
\right)\omega ^5\cosh \left[ {\left( {x-6} \right)\ln \omega }
\right]d\omega } =\xi _\tau \left( 0 \right)\prod\limits_{x_m } {\left(
{1-\frac{x}{x_m }} \right)} .
\end{array}
\end{equation}
And,
\begin{equation}
\label{eq89}
 \begin{aligned}
 \xi _\tau \left( 1 \right)&=2\int\limits_0^1 {L\left( {i\omega }
\right)\omega ^5\cosh \left[ {\left( {1-6} \right)\ln \omega }
\right]d\omega } \\
 &=2\int\limits_0^1 {L\left( {i\omega } \right)\omega ^{10}d\omega } =\xi
_\tau \left( 0 \right)\prod\limits_{x_m } {\left( {1-\frac{5}{x_m }}
\right)} >0,
 \end{aligned}
\end{equation}
which implies that
\begin{equation}
\label{eq90}
\aleph =\sum\limits_{x_m } {\frac{1}{x_m }} <\infty .
\end{equation}
Thus,
\begin{equation}
\label{eq91}
b_0 =\ln \xi _\tau \left( 0 \right)-6\aleph
\end{equation}
is a constant.
\end{remark}

\subsubsection{Two lines of symmetry}

\begin{corollary}
\label{CORO5}
Suppose $\varsigma \in \mathfrak{A}$ and $\varsigma \ne x_m $. Then we have the
following equivalent representations:
\begin{itemize}
  \item
(C1) There is
\begin{equation}
\label{eq92}
\begin{array}{l}
\xi _\tau \left( x \right)=\xi _\tau \left( \varsigma
\right)\prod\limits_{\Im \left( {x_m } \right)>0} {\left[ {1-\frac{\left(
{x-6} \right)^2-\left( {\varsigma -6} \right)^2}{\left( {x_m -6}
\right)^2-\left( {\varsigma -6} \right)^2}} \right]} .
\end{array}
\end{equation}
  \item
(C2) There is
\begin{equation}
\label{eq93}
\begin{array}{l}
\xi _\tau \left( x \right)=\xi _\tau \left( 6 \right)\prod\limits_{\Im
\left( {x_m } \right)>0} {\left[ {1-\left( {\frac{x-6}{x_m -6}} \right)^2}
\right]}.
\end{array}
\end{equation}
\end{itemize}
\end{corollary}

\begin{proof}
By (\ref{eq56}) in Theorem \ref{TH4}, we obtain
\begin{equation}
\label{eq94}
 \begin{aligned}
 \xi _\tau \left( x \right)&=\xi _\tau \left( 0 \right)\prod\limits_{x_m }
{\left( {1-\frac{x}{x_m }} \right)} \\
&=\xi _\tau \left( 0
\right)\prod\limits_{x_m } {\left( {\frac{x_m -x}{x_m }} \right)}\\
&=\xi _\tau
\left( 0 \right)\prod\limits_{x_m } {\left( {\frac{x_m -\varsigma }{x_m
-\varsigma }\cdot \frac{x_m -x}{x_m }} \right)} \\
&=\xi _\tau \left( 0 \right)\prod\limits_{x_m } {\left( {\frac{x_m
-\varsigma }{x_m }\cdot \frac{x_m -x}{x_m -\varsigma }} \right)} \\
&=\xi _\tau
\left( 0 \right)\prod\limits_{x_m } {\frac{x_m -\varsigma }{x_m }}
\prod\limits_{x_m } {\frac{x_m -x}{x_m -\varsigma }}\\
&=\xi _\tau \left( 0
\right)\prod\limits_{x_m } {\frac{x_m -\varsigma }{x_m }} \prod\limits_{x_m
} {\frac{x_m -\varsigma -\left( {x-\varsigma } \right)}{x_m -\varsigma }} \\
&=\xi _\tau \left( 0 \right)\prod\limits_{x_m } {\left( {1-\frac{\varsigma
}{x_m }} \right)} \prod\limits_{x_m } {\frac{x_m -\varsigma -\left(
{x-\varsigma } \right)}{x_m -\varsigma }} \\
&=\xi _\tau \left( 0
\right)\prod\limits_{x_m } {\left( {1-\frac{\varsigma }{x_m }} \right)}
\prod\limits_{x_m } {\left( {1-\frac{x-\varsigma }{x_m -\varsigma }}
\right)} . \\
 \end{aligned}
\end{equation}
Because of
\[
\begin{array}{l}
\xi _\tau \left( \varsigma \right)=\xi _\tau \left( 0
\right)\prod\limits_{x_m } {\left( {1-\frac{\varsigma }{x_m }} \right)} ,
\end{array}
\]
we obtain from (\ref{eq94}) that
\begin{equation}
\label{eq95}
\begin{array}{l}
\xi _\tau \left( x \right)=\xi _\tau \left( \varsigma
\right)\prod\limits_{x_m } {\left( {1-\frac{x-\varsigma }{x_m -\varsigma }}
\right)} .
\end{array}
\end{equation}
Using the functional equation (\ref{eq13}), we reconsider (\ref{eq95}) as
\begin{equation}
\label{eq96}
\begin{array}{l}
 \xi _\tau \left( x \right)=\xi _\tau \left( \varsigma
\right)\prod\limits_{x_m } {\left( {1-\frac{x-\varsigma }{x_m -\varsigma }}
\right)} \\
 =\xi _\tau \left( \varsigma \right)\prod\limits_{\Im \left( {x_m }
\right)>0} {\left( {1-\frac{x-\varsigma }{x_m -\varsigma }} \right)\left(
{1-\frac{x-\varsigma }{12-x_m -\varsigma }} \right)} \\
 =\xi _\tau \left( \varsigma \right)\prod\limits_{\Im \left( {x_m }
\right)>0} {\frac{\left( {x_m -6} \right)-\left( {x-6} \right)}{\left( {x_m
-6} \right)-\left( {\varsigma -6} \right)}} \prod\limits_{\Im \left( {x_m }
\right)>0} {\frac{\left( {6-x_m } \right)-\left( {x-6} \right)}{\left(
{6-x_m } \right)+\left( {6-\varsigma } \right)}} \\
 \end{array}
\end{equation}
where
\begin{equation}
\label{eq97}
 \begin{aligned}
 \prod\limits_{\Im \left( {x_m } \right)>0} {\left( {1-\frac{x-\varsigma
}{x_m -\varsigma }} \right)} &=\prod\limits_{\Im \left( {x_m } \right)>0}
{\left( {1-\frac{\left( {x-6} \right)-\left( {\varsigma -6} \right)}{\left(
{x_m -6} \right)-\left( {\varsigma -6} \right)}} \right)} \\
&=\prod\limits_{\Im \left( {x_m } \right)>0} {\frac{\left[ {\left( {x_m -6}
\right)-\left( {\varsigma -6} \right)} \right]-\left[ {\left( {x-6}
\right)-\left( {\varsigma -6} \right)} \right]}{\left( {x_m -6}
\right)-\left( {\varsigma -6} \right)}} \\
&=\prod\limits_{\Im \left( {x_m } \right)>0} {\frac{\left( {x_m -6}
\right)-\left( {x-6} \right)}{\left( {x_m -6} \right)-\left( {\varsigma -6}
\right)}} \\
 \end{aligned}
\end{equation}
and
\begin{equation}
\label{eq98}
 \begin{aligned}
 \prod\limits_{\Im \left( {x_m } \right)>0} {\left( {1-\frac{x-\varsigma
}{12-x_m -\varsigma }} \right)} &=\prod\limits_{\Im \left( {x_m } \right)>0}
{\left( {1-\frac{\left( {x-6} \right)-\left( {\varsigma -6} \right)}{\left(
{6-x_m } \right)+\left( {6-\varsigma } \right)}} \right)} \\
&=\prod\limits_{\Im \left( {x_m } \right)>0} {\frac{\left[ {\left( {6-x_m }
\right)+\left( {6-\varsigma } \right)} \right]-\left[ {\left( {x-6}
\right)-\left( {\varsigma -6} \right)} \right]}{\left( {6-x_m }
\right)+\left( {6-\varsigma } \right)}} \\
&=\prod\limits_{\Im \left( {x_m } \right)>0} {\frac{\left( {6-x_m }
\right)-\left( {x-6} \right)}{\left( {6-x_m } \right)+\left( {6-\varsigma }
\right)}} . \\
 \end{aligned}
\end{equation}
Further,
\begin{equation}
\label{eq99}
 \begin{aligned}
 \xi _\tau \left( x \right)&=\xi _\tau \left( \varsigma
\right)\prod\limits_{\Im \left( {x_m } \right)>0} {\left[ {\frac{\left( {x_m
-6} \right)-\left( {x-6} \right)}{\left( {x_m -6} \right)-\left( {\varsigma
-6} \right)}\cdot \frac{\left( {6-x_m } \right)-\left( {x-6} \right)}{\left(
{6-x_m } \right)+\left( {6-\varsigma } \right)}} \right]} \\
&=\xi _\tau \left( \varsigma \right)\prod\limits_{\Im \left( {x_m }
\right)>0} {\left[ {\frac{\left( {x_m -6} \right)-\left( {x-6}
\right)}{\left( {x_m -6} \right)-\left( {\varsigma -6} \right)}\cdot
\frac{\left( {x_m -6} \right)+\left( {x-6} \right)}{\left( {x_m -6}
\right)+\left( {\varsigma -6} \right)}} \right]} \\
&=\xi _\tau \left( \varsigma \right)\prod\limits_{\Im \left( {x_m }
\right)>0} {\frac{\left( {x_m -6} \right)^2-\left( {x-6} \right)^2}{\left(
{x_m -6} \right)^2-\left( {\varsigma -6} \right)^2}} . \\
 \end{aligned}
\end{equation}
Finally,
\begin{equation}
\label{eq100}
 \begin{aligned}
 \xi _\tau \left( x \right)&=\xi _\tau \left( \varsigma
\right)\prod\limits_{\Im \left( {x_m } \right)>0} {\frac{\left( {x_m -6}
\right)^2-\left( {x-6} \right)^2}{\left( {x_m -6} \right)^2-\left(
{\varsigma -6} \right)^2}} \\
&=\xi _\tau \left( \varsigma \right)\prod\limits_{\Im \left( {x_m }
\right)>0} {\frac{\left[ {\left( {x_m -6} \right)^2-\left( {\varsigma -6}
\right)^2} \right]-\left[ {\left( {x-6} \right)^2-\left( {\varsigma -6}
\right)^2} \right]}{\left( {x_m -6} \right)^2-\left( {\varsigma -6}
\right)^2}} \\
&=\xi _\tau \left( \varsigma \right)\prod\limits_{\Im \left( {x_m }
\right)>0} {\left[ {1-\frac{\left( {x-6} \right)^2-\left( {\varsigma -6}
\right)^2}{\left( {x_m -6} \right)^2-\left( {\varsigma -6} \right)^2}}
\right]} , \\
 \end{aligned}
\end{equation}
which is the require result.

On the substituting $\varsigma =6$ into (\ref{eq100}), we obtain
\begin{equation}
\label{eq101}
\begin{array}{l}
\xi _\tau \left( x \right)=\xi _\tau \left( \varsigma
\right)\prod\limits_{\Im \left( {x_m } \right)>0} {\left[ {1-\frac{\left(
{x-6} \right)^2}{\left( {x_m -6} \right)^2}} \right]} .
\end{array}
\end{equation}
Thus, the proof is finished.
\end{proof}

\begin{remark}
From (C1) in Corollary \ref{CORO5} we discover that $\xi _\tau \left( x
\right)$ has two lines of symmetry as follows: $x=6$ and $\varsigma =6$. With
(C2) in Corollary \ref{CORO5}, we also discover that we can remove the line
$\varsigma =6$ of symmetry of $\xi _\tau \left( x \right)$ to obtain the
equation (\ref{eq101}) by putting $\varsigma =6$ into (\ref{eq100}).
\end{remark}

As a direct result of Corollary \ref{CORO5}, we also obtain the following:

\begin{corollary}
\label{CORO6}
There is
\begin{equation}
\label{eq2.77a}
\begin{array}{l}
\Lambda \left( \varphi \right)=\Lambda \left( 0 \right)\prod\limits_{\Re
\left( {\varphi _m } \right)>0} {\left( {1-\frac{\varphi ^2}{\varphi _m^2 }}
\right)} ,
\end{array}
\end{equation}
where the product runs over all zeros $\varphi _m $ of $\Lambda \left(
\varphi \right)$.
\end{corollary}

\begin{proof}
Taking $x=6+i\varphi $ into (\ref{eq93}), we have
\begin{equation}
\label{eq2.77b}
\begin{array}{l}
\Lambda \left( \varphi \right)=\xi _\tau \left( {6+i\varphi } \right)=\xi
_\tau \left( 6 \right)\prod\limits_{\Im \left( {x_m } \right)>0} {\left[
{1+\frac{\varphi ^2}{\left( {x_m -6} \right)^2}} \right]}.
\end{array}
\end{equation}

From (\ref{eq2.77b}) all zeros $\varphi _m $ of $\Lambda \left( \varphi \right)$ are
given as follows:
\begin{equation}
\label{eq2.77c}
\begin{array}{l}
x_m -6=\pm i\varphi _m.
\end{array}
\end{equation}

Substituting (\ref{eq2.77c}) into (\ref{eq2.77b}), we show
\begin{equation}
\label{eq2.77d}
\begin{array}{l}
\Lambda \left( \varphi \right)=\xi _\tau \left( 6 \right)\prod\limits_{\Im
\left( {x_m } \right)>0} {\left[ {1+\frac{\varphi ^2}{\left( {x_m -6}
\right)^2}} \right]} =\xi _\tau \left( 6 \right)\prod\limits_{\Re \left(
{\varphi _m } \right)>0} {\left( {1-\frac{\varphi ^2}{\varphi _m^2 }}
\right)}.
\end{array}
\end{equation}

We easily see that the sequence $\left\{ {\varphi _m } \right\}$ is the set of non-zero
zeros of $\Lambda \left( \varphi \right)$ with $|\varphi _{m+1}|>|\varphi _{m}|$.

By (\ref{eq2.77b}), we obtain the relation
\begin{equation}
\label{eq2.77e}
\Lambda \left( \varphi \right)=\xi _\tau \left( {6+i\varphi } \right),
\end{equation}
which leads to

\begin{equation}
\label{eq2.77f}
\Lambda \left( 0 \right)=\xi _\tau \left( 6 \right).
\end{equation}

Putting (\ref{eq2.77f}) into (\ref{eq2.77d}) implies the required result.
\end{proof}

\section{The Riemann-Hardy hypothesis}

\subsection{The conjecture of de Bruijn}

At first, we have the following:

\begin{theorem}
(\textbf{The conjecture of de Bruijn})
\label{TH5}
All of the zeros of $X\left( z \right)$ lie on the critical line $\Re \left(
z \right)=0$.
\end{theorem}

\begin{proof}
By Corollary \ref{CORO1} and Theorem \ref{TH3}, we set up a class of
$X\left( z \right)$ as follows:
\begin{equation}
\label{eq102}
\begin{array}{l}
X\left( z \right)=\sum\limits_{s=0}^\infty {\chi \left( s \right)z^{2s}}
=X\left( 0 \right)\prod\limits_{\Im \left( {z_m } \right)>0} {\left(
{1-\frac{z^2}{z_m^2 }} \right)} .
\end{array}
\end{equation}
From (\ref{eq28}) and (\ref{eq102}) we have
\begin{equation}
\label{eq103}
X\left( z \right)=\sum\limits_{s=0}^\infty {\chi \left( s \right)z^{2s}}
\end{equation}
such that
\begin{equation}
\label{eq104}
\overline {X\left( z \right)} =\overline {\left[ {\sum\limits_{s=0}^\infty
{\chi \left( s \right)z^{2s}} } \right]} =\sum\limits_{s=0}^\infty {\chi
\left( s \right)\overline z ^{2s}} .
\end{equation}
in which $\overline {X\left( z \right)} $ and $\overline z $ are the complex
conjugates of $X\left( z \right)$ and $z$, respectively.

It follows from (\ref{eq104}) that
\begin{equation}
\label{eq105}
\overline {X\left( z \right)} =X\left( {\overline z } \right).
\end{equation}
On adopting (\ref{eq102}) and (\ref{eq105}), we obtain the first product of $X\left(
{\overline z } \right)$ as follows:
\begin{equation}
\label{eq106}
\begin{array}{l}
\overline {X\left( z \right)} =X\left( {\overline z } \right)=X\left( 0
\right)\prod\limits_{\Im \left( {z_m } \right)>0} {\left( {1-\frac{\overline
z ^2}{z_m^2 }} \right)} .
\end{array}
\end{equation}
Similarly, by finding the complex conjugate of $X\left( z \right)$in (\ref{eq102})
and using (\ref{eq47}), we get
\begin{equation}
\label{eq107}
 \begin{aligned}
 \overline {X\left( z \right)} &=\overline {\left[ {X\left( 0
\right)\prod\limits_{\Im \left( {z_m } \right)>0} {\left(
{1-\frac{z^2}{z_m^2 }} \right)} } \right]}
&=X\left( 0 \right)\overline
{\left[ {\prod\limits_{\Im \left( {z_m } \right)>0} {\left(
{1-\frac{z^2}{z_m^2 }} \right)} } \right]} \\
&=X\left( 0 \right)\prod\limits_{\Im \left( {z_m } \right)>0} {\left(
{1-\frac{\overline z ^2}{\overline {z_m } ^2}} \right)} . \\
 \end{aligned}
\end{equation}
Combining (\ref{eq104}) and (\ref{eq107}), we present
\begin{equation}
\label{eq108}
\begin{array}{l}
\overline {X\left( z \right)} =X\left( {\overline z } \right)=X\left( 0
\right)\prod\limits_{\Im \left( {z_m } \right)>0} {\left( {1-\frac{\overline
z ^2}{z_m^2 }} \right)} =X\left( 0 \right)\prod\limits_{\Im \left( {z_m }
\right)>0} {\left( {1-\frac{\overline z ^2}{\overline {z_m } ^2}} \right)}.
\end{array}
\end{equation}
Taking the complex conjugate of (\ref{eq108}), we may get
\begin{equation}
\label{eq109}
\begin{array}{l}
X\left( z \right)=X\left( 0 \right)\prod\limits_{\Im \left( {z_m }
\right)>0} {\left( {1-\frac{z^2}{z_m^2 }} \right)} =X\left( 0
\right)\prod\limits_{\Im \left( {z_m } \right)>0} {\left(
{1-\frac{z^2}{\overline {z_m } ^2}} \right)} .
\end{array}
\end{equation}
By Corollary \ref{CORO1}, we know that there exists any positive number
$\varepsilon >0$ such that the series
\begin{equation}
\label{eq110}
\sum\limits_{m=1}^\infty {\left| {z_m } \right|^{-\left( {1+\varepsilon }
\right)}}
\end{equation}
is convergent.

This implies that
\begin{equation}
\label{eq111}
\sum\limits_{m=1}^\infty {\left| {\overline {z_m } } \right|^{-\left(
{1+\varepsilon } \right)}}
\end{equation}
is convergent.

From (\ref{eq110}) and (\ref{eq111}),
\begin{equation}
\label{eq112}
\sum\limits_{m=1}^\infty {\left| {z_m } \right|^{-2}}
\end{equation}
and
\begin{equation}
\label{eq113}
\sum\limits_{m=1}^\infty {\left| {\overline {z_m } } \right|^{-2}}
\end{equation}
are convergent.

Putting $z=1$ into (\ref{eq109}), we present
\begin{equation}
\label{eq114}
\begin{array}{l}
X\left( 1 \right)=X\left( 0 \right)\prod\limits_{\Im \left( {z_m }
\right)>0} {\left( {1-\frac{1}{z_m^2 }} \right)} =X\left( 0
\right)\prod\limits_{\Im \left( {z_m } \right)>0} {\left(
{1-\frac{1}{\overline {z_m } ^2}} \right)} .
\end{array}
\end{equation}

By using the fact (\ref{eq112}) and (\ref{eq113}) are convergent, Theorem 5 in Knopp's
monograph (see \cite{29}, p.10) said that (\ref{eq114}) is absolutely convergent.
Following Theorem 3 in Knopp's monograph (see \cite{29}, p.10), we see that
(\ref{eq114}) converges.

Making use of Theorem 3 in Knopp's monograph (see \cite{29}, p.10) and (\ref{eq114}),
this implies that
\begin{equation}
\label{eq115}
\sum\limits_{m=1}^\infty {\frac{1}{z_m^2 }}
\end{equation}
and
\begin{equation}
\label{eq116}
\sum\limits_{m=1}^\infty {\frac{1}{\overline {z_m } ^2}}
\end{equation}
are convergent, and
\begin{equation}
\label{eq117}
\sum\limits_{m=1}^\infty {\frac{1}{z_m^2 }} =\sum\limits_{m=1}^\infty
{\frac{1}{\overline {z_m } ^2}} .
\end{equation}
From (\ref{eq117}) we have
\begin{equation}
\label{eq118}
z_m^2 -\overline {z_m } ^2=0
\end{equation}
such that
\begin{equation}
\label{eq119}
\left( {z_m -\overline {z_m } } \right)\left( {z_m +\overline {z_m } }
\right)=0.
\end{equation}
By using the fact $z_m -\overline {z_m } =2\Im \left( {z_m } \right)\ne 0$,
it follows from (\ref{eq119}) that
\begin{equation}
\label{eq120}
z_m +\overline {z_m } =2\Re \left( {z_m } \right)=0.
\end{equation}
Thus,
\begin{equation}
\label{eq121}
\Re \left( {z_m } \right)=0
\end{equation}
implies that
\begin{equation}
\label{eq122}
z_m =i\Im \left( {z_m } \right).
\end{equation}
Substituting (\ref{eq122}) back into (\ref{eq117}) leads to the series
\[
\sum\limits_{m=1}^\infty {\frac{1}{z_m^2 }} =\sum\limits_{m=1}^\infty
{\frac{1}{\overline {z_m } ^2}} =-\sum\limits_{m=1}^\infty {\frac{1}{\Im
\left( {z_m } \right)^2}}
\]
converges.

Applying Theorem 4 in Knopp's monograph (see \cite{29}, p.10), this implies that
(\ref{eq114}) converge. Once again, it implies that we go back to verify the truth
of (\ref{eq114}) and (\ref{eq122}).

Let $\left| {\Im \left( {z_m } \right)} \right|=\gamma _m >0$. Then we have
from Theorem \ref{TH3} that
\begin{equation}
\label{eq123}
\begin{array}{l}
X\left( z \right)=X\left( 0 \right)\prod\limits_{\Im \left( {z_m }
\right)>0} {\left( {1-\frac{z^2}{z_m^2 }} \right)} =X\left( 0
\right)\prod\limits_{m=1}^\infty {\left( {1+\frac{z^2}{\gamma _m^2 }}
\right)} .
\end{array}
\end{equation}
Thus, we complete the proof.
\end{proof}

\subsection{The conjecture of Ki}

We now present the following result:

\begin{theorem}
(\textbf{The conjecture of Ki})
\label{TH6}
All of the zeros of $\Lambda \left( \varphi \right)$ lie on the critical
line $\Im \left( \varphi \right)=0$.
\end{theorem}

\begin{proof}
By using (\ref{eq123}), we have
\begin{equation}
\label{eq124}
\Lambda \left( \varphi \right)=X\left( {i\varphi } \right)
\end{equation}
such that
\begin{equation}
\label{eq125}
\begin{array}{l}
\Lambda \left( \varphi \right)=X\left( {i\varphi } \right)=X\left( 0
\right)\prod\limits_{m=1}^\infty {\left( {1-\frac{\varphi ^2}{\gamma _m^2 }}
\right)} .
\end{array}
\end{equation}
By using $\Lambda \left( 0 \right)=X\left( 0 \right)$, (\ref{eq125}) can be
rewritten as
\begin{equation}
\label{eq126}
\begin{array}{l}
\Lambda \left( \varphi \right)=\Lambda \left( 0
\right)\prod\limits_{m=1}^\infty {\left( {1-\frac{\varphi ^2}{\gamma _m^2 }}
\right)} .
\end{array}
\end{equation}

From (\ref{eq125}) we obtain the required result.

This proof of Theorem \ref{TH6} is now finished.
\end{proof}

\subsection{The proof of Theorem \ref{TH1}}
We now return to its proof.

Applying (\ref{eq54}) and (\ref{eq122}), we obtain
\begin{equation}
\label{eq127}
\begin{array}{l}
\xi _\tau \left( x \right)=\xi _\tau \left( 6 \right)\prod\limits_{\Im
\left( {z_m } \right)>0} {\left[ {1-\frac{\left( {x-6} \right)^2}{z_m^2 }}
\right]} =\xi _\tau \left( 6 \right)\prod\limits_{\Im \left( {z_m }
\right)>0} {\left[ {1+\frac{\left( {x-6} \right)^2}{\left| {\Im \left( {z_m
} \right)} \right|^2}} \right]} .
\end{array}
\end{equation}
On substituting $\left| {\Im \left( {z_m } \right)} \right|=\gamma _m >0$
into (\ref{eq127}), this gives
\begin{equation}
\label{eq128}
\begin{array}{l}
\xi _\tau \left( x \right)=\xi _\tau \left( 6 \right)\prod\limits_{\Im
\left( {z_m } \right)>0} {\left[ {1+\frac{\left( {x-6} \right)^2}{\left|
{\Im \left( {z_m } \right)} \right|^2}} \right]} =\xi _\tau \left( 6
\right)\sum\limits_{m=1}^\infty {\left[ {1+\frac{\left( {x-6}
\right)^2}{\gamma _m^2 }} \right]} .
\end{array}
\end{equation}
Adopting Theorem 1 in Knopp's book (see \cite{29}, p.9) and taking $\xi _\tau
\left( x \right)=0$ in (\ref{eq128}),
\begin{equation}
\label{eq129}
1+\frac{\left( {x-6} \right)^2}{\gamma _m^2 }=0,
\end{equation}
is the root of $\xi _\tau \left( x \right)=0$. Thus,
\begin{equation}
\label{eq131}
x_m =6\pm i\gamma _m ,
\end{equation}
where $\gamma _m >0$.

From (\ref{eq7}) and (\ref{eq131}), we deduce that all zeros of $\xi _\tau \left( x
\right)$ lies on $\Re \left( {x_m } \right)=6$.

Thus, we complete the proof of Theorem \ref{TH1}.


\begin{thebibliography}{120}

\bibitem{1}
S. Ramanujan, On certain arithmetical functions, Transactions of the Cambridge Philosophical Society, 22 (1916) (9), 159-184.

\bibitem{2}
G. H. Hardy, Note on Ramanujan's arithmetical function $\tau (n)$, Mathematical Proceedings of the Cambridge Philosophical Society, 23 (1927) (06), 675.

\bibitem{3}
L. J. Mordell, On Mr Ramanujan's empirical expansions of modular functions, Proceedings of the Cambridge Philosophical Society, 19(1917), 117-124.

\bibitem{4}
W. C. Winnie Li, The Ramanujan conjecture and its applications, Philosophical Transactions of the Royal Society A, 378 (2020)(2163), 20180441.

\bibitem{5}
R. A. Rankin, Contributions to the theory of Ramanujan's function $\tau(n)$ and similar arithmetical functions: Ii. the order of the Fourier coefficients of integral modular forms, Mathematical Proceedings of the Cambridge Philosophical Society, 35(1939) (3), pp. 357-372.

\bibitem{6}
M. Rogers, Identities for the Ramanujan zeta function, Advances in Applied Mathematics, 51(2013) (2), 266-275.

\bibitem{7}
B. Heim, M. Neuhauser and F. Rupp, Fourier coefficients of powers of the Dedekind eta function, The Ramanujan Journal, 48 (2019) (1), 1-11.

\bibitem{8}
P. Deligne, La conjecture de Weil. I., Publications Math\'{e}matiques de l'Institut des Hautes \'{E}tudes Scientifiques, 43 (1974)(1), 273-307.

\bibitem{9}
J. S. Balakrishnan, W. Craig and K. Ono, Variations of Lehmer's conjecture for Ramanujan's tau-function, Journal of Number Theory, 237 (2022), 3-14.

\bibitem{10}
C. J. Moreno, Prime number theorems for the coefficients of modular forms, Bulletin of the American Mathematical Society, 78 (1972) (5), 796-798.

\bibitem{11}
R. A. Rankin, Contributions to the theory of Ramanujan's function $\tau (n)$ and similar arithmetical functions: I. The zeros of the function on the line, Mathematical Proceedings of the Cambridge Philosophical Society, 35 (1939) (3), 351-356.

\bibitem{12}
J. B. Conrey and A. Ghosh, Tur\'{a}n inequalities and zeros of Dirichlet series associated with certain cusp forms, Transactions of the American Mathematical Society, 342 (1994) (1), 407-419.

\bibitem{13}
J. R. Wilton, A note on Ramanujan's arithmetical function $\tau (n)$, Mathematical Proceedings of the Cambridge Philosophical Society, 25(1929) (2), 121-129.

\bibitem{13a}
B. C. Berndt and M. I. Knopp, Hecke's theory of modular forms and Dirichlet series, Vol. 5, World Scientific, 2008.

\bibitem{14}
J. B. Conrey and A. Ghosh, Simple zeros of the Ramanujan $\tau $-Dirichlet series, Inventiones mathematicae, 94 (1988) (2), 403-419.

\bibitem{15}
N. G. De Bruijn, The roots of trigonometric integrals, Duke Mathematical Journal, 17 (1950) (3), 197-226.

\bibitem{16}
G. H. Hardy, Ramanujan twelve lectures on subjects suggested by his life and work, Cambridge University Press, Cambridge, 1940.

\bibitem{17}
R. Courant and F. John, Introduction to Calculus and Analysis I, Springer, 2012.

\bibitem{18}
C. G. Lekkerkerker, On the zeros of a class of Dirichlet series, van Gorcum, 1955.

\bibitem{19}
B. C. Berndt, On the zeros of a class of Dirichlet series I, Illinois Journal of Mathematics, 14 (1970) (2), 244-258.

\bibitem{20}
J. L. Hafner, On the zeros of Dirichlet series associated with certain cusp forms, Bulletin of the American Mathematical Society, 8 (1983) (2), 340-342.

\bibitem{21}
H. Ki, On the zeros of approximations of the Ramanujan $\Xi $-function, The Ramanujan Journal, 17 (2008) (1), 123-143.

\bibitem{22}
A. Chirre and O.V. Casta\~{n}\'{o}n, A note on the zeros of approximations of the Ramanujan $\Xi $- function, The Ramanujan Journal, 57(2020), 389-400.

\bibitem{23}
H. R. P. Ferguson, R. D. Major, K. E. Powell and H. G. Throolin, On zeros of mellin transforms of
$SL_2(\mathbf{Z})$ cusp forms, Mathematics of Computation, 42 (1984) (165), 241.

\bibitem{24}
J. B. Keiper, On the zeros of the Ramanujan $\tau$-Dirichlet series in the critical strip, Mathematics of Computation, 65 (1996) (216), 1613-1620.

\bibitem{25}
P. Sarnak, Some applications of modular forms, Vol. 99, Cambridge University Press, 1990.

\bibitem{26}
T. M. Apostol and A. Sklar, The approximate functional equation of Hecke's Dirichlet series, Transactions of the American Mathematical Society, 86 (1957) (2), 446-462.

\bibitem{26a}
A. Ogg, Modular forms and Dirichlet series, Vol. 39, WA Benjamin, New York, 1969.

\bibitem{27}
B. Y. Levin, Distribution of zeros of entire functions, Vol. 150, American Mathematical Society, 1980.

\bibitem{28}
E. C. Titchmarsh, The theory of functions, Oxford University Press, 1939.

\bibitem{29}
K. Knopp, Theory of functions, Parts II, Dover Publications, New York, 1947.
\end{thebibliography}
\end{document}